\documentstyle[amssymb]{amsart}

\def\beq{\begin{equation}}
\def\eeq{\end{equation}}

\def\bq{\begin{quote}}
\def\eq{\end{quote}}

\def\lqq{\lq \lq }
\relax

\theoremstyle{definition}

\newcommand{\be}{\begin{equation}}
\newcommand{\ee}{\end{equation}}
\begin{document}

\title[On some operator identities \dots]
{On some operator identities and representations of algebras}
\author{ Alexander Turbiner and Gerhard Post}
\address
{Institute for Theoretical and Experimental Physics,
Moscow 117259, Russia \quad {\it and}
Deparment of Applied Mathematics, University of Twente, P.O. Box 217,
\quad \quad 7553 CV Enschede, Holland}
\email{ TURBINER@@CERNVM or TURBINER@@VXCERN.CERN.CH \quad \quad
{\it and}
 POST@@MATH.UTWENTE.NL}
\curraddr{I.H.E.S., Bures-sur-Yvette, F-91440, France}

\date{}
\maketitle
\begin{abstract}
Certain infinite families of operator identities related to powers of
positive root generators of (super) Lie algebras of first-order
differential operators and $q$-deformed algebras of first-order
finite-difference operators are presented.
\end{abstract}

1. The following {\it operator identity} holds
\label{e1}
\begin{equation}
{(J^+_n)}^{n+1} \equiv (x^2 \partial_x - n x)^{n+1} =
x^{2n+2}\partial_x^{n+1} , \partial_x \equiv {d \over dx}, n=0,1,2,\ldots
\end{equation}
The proof is straightforward:
\begin{enumerate}
\item[(i)] the operator ${(J^+_n)}^{n+1}$
annihilates the space of all polynomials of degree not higher than $n$,
${\cal P}_n(x)=Span\{x^i: 0 \leq i \leq n\}$;
\item[(ii)] in general, an $(n+1)-$th order linear differential
operator annihilating ${\cal P}_n(x)$ must have the form
$B(x)\partial_x^{n+1}$, where $B(x)$ is an arbitrary function and
\item[(iii)]
since ${(J^+_n)}^{n+1}$ is a graded operator, deg$(J^+_n)=+1$,
/footnote{so $J^+_n$ maps $x^k$ to a multiple of $x^{k+1}$ }
deg$({J^+_n}^{n+1})=n+1$, hence $B(x)=b x^{2n+2}$ while clearly the
constant $b=1$.
\end{enumerate}
It is worth noting that taking the degree in (1) different from
 $(n+1)$, the l.h.s. in (1) will contain all derivative terms from zero up
to $(n+1)$-th order.

The identity has a Lie-algebraic interpretation.
The operator $(J^+_n)$ is the positive-root generator of the algebra $sl_2$
of first-order differential operators (the other $sl_2$-generators are
$J^0_n = x \partial_x - n/2 \ , J^-_n = \partial_x $). Correspondingly, the
space ${\cal P}_n(x)$ is nothing but the $(n+1)$-dimensional
irreducible representation of $sl_2$. The identity (1) is a
consequence of the fact that ${(J^+_n)}^{n+1}=0$ in this
representation.

There exist other algebras of differential or finite-difference
operators (in more than one variable), which admit a
finite-dimensional representation.  This leads to more general
and remarkable operator identities. In the present Note, we show
that (1) is one representative of infinite
family of identities for differential and finite-difference
operators.


2. The Lie-algebraic interpretation presented above allows us to generalize
(1) for the case of differential operators of several variables,
taking appropriate degrees of the highest-positive-root generators of
(super) Lie algebras of first-order differential operators, possessing
a finite-dimensional invariant sub-space (see e.g.\cite{t1}).
First we consider
the case of $sl_3$. There exists a representation of $sl_3(\bold C)$ as
differential operators on $\bold C^2$. One of the generators is
\[ J^1_2 (n)= x^2 \partial_x\ +\ xy \partial_y - n x \]
The space ${\cal P}_n(x,y)=Span\{x^iy^j: 0 \leq i+j \leq n\}$ is a
finite-dimensional representation for $sl_3$, and due to the fact $(J^1_2
(n))^{n+1}=0$ on the space ${\cal P}_n(x,y)$, hence we have
\label{e2}
\begin{eqnarray}
{(J^1_2 (n))}^{n+1} = (x^2 \partial_x\ +\ xy \partial_y - n x)^{n+1} =
\nonumber \\
\sum_{k=0}^{k=n+1} {n+1 \choose k} x^{2n+2-k}y^k \partial_x^{n+1-k}
\partial_y^k  \ ,
\end{eqnarray}
This identity is valid for $y \in \bold C$ (as described above), but also if
$y$ is a Grassmann variable, i.e. $y^2=0$ \footnote{In this case just
two terms in l.h.s. of (2) survive.}. In the last case, $J^1_2 (n)$ is
a generator of $osp(2,2)$, see \cite{t1}.

More general (using $sl_k$ instead of $sl_3$), the following
operator identity holds
\label{e3}
\begin{eqnarray}
{(J^{k-2}_{k-1} (n))}^{n+1} \equiv (x_1\sum_{m=1}^{k} (x_m \partial_{x_m}\
- n))^{n+1} =
\nonumber \\
x_1^{n+1}\sum_{j_1+j_2+\ldots+j_k=n+1} C^{n+1}_{j_1,j_2,\ldots,j_k}
x_1^{j_1}x_2^{j_2}\ldots x_k^{j_k} \partial_{x_1}^{j_1}
\partial_{x_2}^{j_2}\ldots \partial_{x_k}^{j_k}   \ ,
\end{eqnarray}
where $C^{n+1}_{j_1,j_2,\ldots,j_k}$ are the generalized binomial
(multinomial) coefficients.
If $x \in \bold C^k$, then $J^{k-2}_{k-1} (n)$ is a generator of the algebra
$sl_k(\bold C)$ \cite{t1}, while some of the variables $x'$s are Grassmann
ones, the operator $J^{k-2}_{k-1} (n)$ is a generator of a certain super
Lie algebra of first-order differential operators. The operator in l.h.s.
of (3) annihilates the linear space of polynomials
${\cal P}_n(x_1,x_2,\ldots x_k)=Span\{x_1^{j_1}x_2^{j_2}\ldots x_k^{j_k} :
 0 \leq j_1+j_2+\ldots+j_k \leq n\}$.


3. The above-described family of operator identities can be generalized
for the case of finite-difference operators with the Jackson symbol, $D_x$
(see e.g. \cite{e})
\[ D_x f(x) = {{f(x) - f(q^2x)} \over {(1 - q^2) x}} + f(q^2x) D_x\]
instead of the ordinary derivative. Here, $q$ is an arbitrary complex
number. The following operator identity holds
\label{e4}
\begin{equation}
{(\tilde J^+_n)}^{n+1} \equiv  (x^2 D_x - \{ n \}  x)^{n+1} =
q^{2n(n+1)}x^{2n+2} D^{n+1}_x , n=0,1,2,\ldots
\end{equation}
(cf.(1)), where $\{n\} = {{1 - q^{2n}}\over {1 - q^2}}$ is so-called
$q$-number.  The operator in the r.h.s. annihilates the space ${\cal
P}_n(x)$. The proof is similar to the proof of the identity (1).

{}From algebraic point of view the operator $\tilde J^+_n$ is the generator
of a $q$-deformed algebra $sl_2(\bold C)_q$ of first-order
finite-difference operators on the line: \linebreak
$\tilde J^0_n = \ x D - \hat{n},\ \tilde J^-_n = \ D $, where
$\hat n \equiv {\{n\}\{n+1\}\over \{2n+2\}}$ (see \cite{ot} and
also \cite{t1}), obeying the commutation relations
\label{e5}
\[
q^2 \tilde  j^0\tilde  j^- \ - \ \tilde  j^-\tilde  j^0 \
= \ - \tilde  j^-
\]
\begin{equation}
 q^4 \tilde  j^+\tilde  j^- \ - \ \tilde  j^-\tilde  j^+ \
= \ - (q^2+1) \tilde  j^0
\end{equation}
\[
\tilde  j^0\tilde  j^+ \ - \ q^2\tilde  j^+\tilde  j^0 \ = \  \tilde  j^+
 \]
($\tilde j$'s are related with $\tilde J$'s through some multiplicative
factors). The algebra (5) has the linear space ${\cal P}_n(x)$ as a
finite-dimensional representation.

An attempt to generalize (2) replacing continuous derivatives by Jackson
symbols immediately leads to requirement tp introduce the quantum plane and
$q$-differential calculus \cite{wz}
\label{e6}
\[ xy=qyx\ , \]
\[ D_x x=1+q^2 xD_x+(q^2-1) yD_y \quad ,\quad D_x y=qyD_x \ ,\]
\[ D_y x=qxD_y\quad , \quad D_y y=1+q^2 yD_y \ ,\]
\begin{equation}
 D_xD_y=q^{-1}D_yD_x \ .
\end{equation}
The formulae analogous to (2) have the form
\label{e7}
\[ {(\tilde J^1_2 (n))}^{n+1} \equiv
(x^2 D_x\ +\ xy D_y - \{ n\} x)^{n+1} = \]
\begin{equation}
\sum_{k=0}^{k=n+1} q^{2n^2-n(k-2)+k(k-1)} {n+1 \choose k}_q x^{2n+2-k}y^k
D_x^{n+1-k} D_y^k  \ ,
\end{equation}
where
\[ {n \choose k}_q \equiv {\{n\}! \over {\{k\}!\{n-k\}!}}\ ,\ \{n\}! =
\{1\} \{2\}\ldots \{n\} \]
are $q$-binomial coefficient and $q$-factorial, respectively. Like all
previous cases, if $y \in {\bold C}$, the operator $\tilde J^1_2 (n)$ is
one of generators of
$q$-deformed algebra $sl_3(\bold C)_q$ of finite-difference
operators, acting on the quantum plane and having the linear space
${\cal P}_n(x,y)=Span\{x^iy^j: 0 \leq i+j \leq n\}$
as a finite-dimensional representation; the l.h.s. of (7) annihilates
${\cal P}_n(x,y)$. If $y$ is Grassmann variable,  $\tilde J^1_2 (n)$ is a
generator of the $q$-deformed superalgebra $osp(2,2)_q$ possessing
finite-dimensional representation.

Introducing a quantum hyperplane \cite{wz}, one can generalize
the whole family of the operator identities (3) replacing continuous
derivatives by finite-difference operators.

One of us (A.T.) acknowledges of Profs. L. Michel, R. Thom
and F. Pham for kind hospitality and their interest to the present work.

\end{document}